\documentclass[12pt]{amsart}
\usepackage{amscd}
\usepackage[arrow,matrix]{xy}
\usepackage{graphicx}
\usepackage{amsmath}
\begin{document}
\title[Application of  integral  geometry to minimal surfaces]{Application of integral geometry to minimal surfaces}
\author{ L\^E H\^ong V\^an}

\abstract This  is a corrected version  of  my    paper ``Application  of integral  geometry  to minimal surfaces"   appeared in  International J. Math. vol. 4 Nr. 1 (1993),
89-111.  The correction  concerns  Proposition 3.5.  We  discuss  this correction in Appendix   to the original  version     of my published paper by reproducing   our correspondence  with Professor  Tasaki.
\endabstract
\maketitle

{\small  CONTENTS\\
1. Introduction.\\
2. General construction and examples. \\
3. Minimal cycles in Grassmannian manifolds.\\
4. Properties of $(M)^*$- minimal cycles.\\
5.  Appendix: Correspondence  with Professor  Tasaki  on Proposition 3.5  \\
References.}

\vspace{1cm}

\S 1. {\bf Introduction}.

The theory of higher dimensional minimal surfaces,
especially its main branch - the Plateau problem, has been intensively developed
since the sixties when E. R. Reifenberg, H. Federer, W. H. Fleming, E. De Giorgi and 
F. Almgren proved existence and almost regularity
theorems for solutions of the higher dimensional Plateau problem (or simply
speaking, globally minimal surfaces) in different contexts of geometric measure
theory. After that, the other part of the theory, namely, construction, classification and 
study of geometry of globally minimal surfaces has been developed rapidly. 
The first non-trivial example of globally
minimal surfaces was obtained by H. Federer by showing that every K\"ahler
submanifold is  globally (homologically) minimal  in its ambient K\"ahler 
manifold [ Fe 1]. His method of employing  exterior powers of the K\"ahler 
form in K\"ahler manifolds has been generalized for other Riemannian manifolds 
in the works of M. Berger, H. B. Lawson, Dao Trong Thi, R. Harvey and H. B. Lawson
([Be], [Ln 2], [D], [H-L]). Now, this method is called 
the calibration method and it has  various applications in the
 study of geometry of  globally minimal surfaces as well as of (locally) minimal surfaces ([DGGW 1],
[DGGW 2], [G-M-Z], [Le 1], [Le 2], [Lr],...).
Other interesting examples of globally minimal surfaces were obtained by A.T. Fomenko
[Fo 1, Le-Fo] by using an estimate from below for the volume of globally minimal surfaces in
Riemannian manifolds. His idea came from  Griffiths' idea of using exhaustion
functions on algebraic manifolds in the Nevalinna theory. His method allows us to construct homological
minimal submanifolds when the coefficient group of homologies may be finite
$({\mathbf Z_p})$ or infinite (${\mathbf Z}$). Note that the calibration method works only for
homology groups  with coefficients in ${\mathbf R}$. But Fomenko's method which
 depends on an estimate
involving only the injective radius, Riemannian curvature of ambient manifolds and dimension of submanifolds, cannot give us so many examples of globally minimal surfaces. To our knowledge, 
up to now, all non-trivial examples of globally minimal surfaces are 
obtained by using the above mentioned methods with the exception of some globally minimal 
hypersurfaces with  a large symmetry groups where one can reduce the problem
of higher dimension to dimension 2 which can be completely analysed. This
reduction method was invented by W-Y. Hsiang and H. B. Lawson [Hs-Ln] and [Ln 1]. 

This paper is an attempt to fill the gap between the calibration method and the Fomenko
method. This new method may be also called an analog of the calibration method for
discrete coefficients of homology groups (of Riemannian manifolds). The idea is
simple; it also comes from complex geometry. Let us recall the Crofton-type
formula (which has originated in probability theory [Sa]).

{\bf Theorem.} [Ch, p.146] {\it Let $f: M \longrightarrow {\mathbf C}P^n $ be a compact holomorphic
curve with or without boundary. Then

$$ \int_{{\mathbf C}P^n} \# (f(M) \bigcap \gamma) \, d \gamma  = Area (M) , \eqno(1.1)$$

where $\gamma$ is a (complex) hyperplane of ${\mathbf C}P^n$, and the space of these hyperplanes
is identified with ${\mathbf C}P^n$ equipped with the invariant measure, and $\# (X)$ denotes
the number of points in $X$. }

A more detailed analysis shows that if we replace a holomorphic curve $M$ by any
(real) two-dimensional surface $M'$, then the equality (1.1) becomes   an inequality, 
where the right hand side is greater than the left one (see Proposition 2.11 and Proposition 3.10 which we call {\it 
Integral Wirtinger Inequality}). So, this strengthened Crofton-type formula
gives us a new proof of homological minimality of ${\mathbf C}P^1$ , and moreover, an estimate on the
measure of all (complex) hyperplanes meeting a fixed holomorphic curve k times
(see Equidistribution Theorem [Ch, p.146] and Theorem 4.1). In fact, some
authors have used similar integral formulae in  order to estimate the volume of
2-dimensional analytical sets in ${\mathbf C}^n$, but their formulae concern only
the simplest case of real dimension 1 (cf.[K-R] and references in that paper). Our idea is a natural generalization of the Crofton-type formula.
Namely, we want to estimate the volume of a submanifold $ N
\subset M$ by its intersection number $ \# (N \bigcap N^* _{\lambda})$, where  
$N^*_{\lambda}$ is a family of submanifolds in $M$. Since the
algebraic intersection number is a homology invariant we hope to get an estimate
from below for the volume of a submanifold realizing a given cycle. The use of
intersection number as a homology invariant explains the analogy between this
method and the calibration method, which essentially employs another homology
invariant - the Stokes formula. But in view of the Federer stability theorem
[Fe 2] the relation between these methods proves to be more intimate; in many cases,
the effectiveness of one method leads to the effectiveness of the other one (see
\S 4). Applying this intersection method we obtain some old and new examples of
globally minimal submanifolds in Grassmannian spaces. In a few cases this gives us
a classification theorem for globally minimal submanifolds in a certain class (see
\S 3 and \S 4) and their new properties  such as equidistribution in
measure of globally minimal surfaces. Other applications of integral geometry
to minimal surfaces will appear in our next paper. The present note is based on a revised form of author's preprint [Le 3].

\vspace{0,5cm}

{\bf \S 2. General construction and examples.}

Let us begin with a simple example. 

{\it Example 2.1}. Let $M^m$ be a Riemannian manifold and $TM$ its tangent bundle.
Let the Riemanian metric on $M$ be naturally lifted on $TM$. Then $M^m$ realizes
a nontrivial cycle in the homology group $H_m (TM, {\mathbf Z_2})$ and moreover it has the
minimal volume in its homology class $[M]$. In fact, if $M'$ is another submanifold in $TM$ and realizing the
cycle $[M] \in H_* (TM, {\mathbf Z_2})$, then $M'$ must meet every fiber $\pi 
 _x,
x \in M$. Consequently, the projection $\pi : M' \longrightarrow M$ is surjective.
It is easy to see that the projection $\pi$ decreases the volume element (in any
dimension not exceeding $\dim M = m $). Hence we get the assertion. This example is
interesting
because if $M$ is not orientable then $H_m (TM, {\mathbf Z}) = 0$ and the classical
calibration method is not applicable!

Now let us give a {\it general construction}, which generally does not depend
on fibrations (such simple fibrations as the above example  occur very
rarely). Let us consider a Riemannian manifold $M^m$. Suppose we have a family
$(M)^*$ of $n$-dimensional submanifolds $N_y \subset M, \; y\in (M)^*$ . Suppose
further that $(M)^*$ is a smooth manifold with a volume element  $\mu_y$= $\text{  vol}_{m^*}$, where $m^*$ is the dimension of $(M)^*$. For
every $ X \subset M$ denote by $S_X \subset (M)^*$ the set of all submanifolds
$N_y$ passing through the set $X$. Now we fix a point $x \in M$ and a $(m-n)$-dimensional
subspace
$V^{m-n} \subset T_x M$. Denote by $B (x,V^{m-n},r)$ the geodesic ball of radius $r$ in
$M$ with its center at $x$ and its tangent space at $x$ equal to $ V^{m-n}$. Let us consider the following limit
$$ cd(x,V^{m-n}) = \lim_{ r \to 0} { \text{  vol} _{m^*}(S_{B (x,V^{m-n},r)}) \over    \text{  vol}
           (B (x,V^{m-n},r)) }  .\eqno(2.1) $$

 Suppose for every $x \in M$ the set $S_x $ is a compact smooth submanifold
in $(M)^*$. Then the limit in (2.1) exists. To compute this
limit we fix a submanifold $S_x$ and a small normal neighborhood of
$S_x$ in $(M)^*$. Obviously, there exists a fiber bundle $F$ over $S_x$ in this neighborhood such that $S_x$ is embedded into it as a zero section of generic position. For
instance, in order to construct $F$ we can use the exponential map from the normal bundle over $S_x$
to $(M)^*$. For every $y\in S_x$, with the help of $F$, we can construct a map $F_y$
from a neighborhood of $x \in M$ to the fixed neighborhood of $S_x$ as follows: $M \ni x' \mapsto S_{x'} \cap p^{-1}y$, where $p^{-1}y$ is the fiber over $y\in S_x$. Since
$S_x$ meets fibers transversally, the map $F_y$ is well defined in a sufficiently small neighborhood of $x$, that is, $p^{-1}y$ meets $S_{x'}$  only at
one point. Then we have
 $$cd(x,V^{m-n}) = \int_{S_x} \text{  vol}(\overline {T_y S_x } \wedge dF_y(\overline{V^{m-n}})).$$
Here for any linear subspace $L$ we denote by $\overline{L}$ the unit polyvector associated with $L$.
We call the limit in (2.1) { \it a deformation coefficient  $cd(x,V^{m-n})$}. Put
$$\sigma (M)^*_{m-n} = \max  \{cd(x,V^{m-n}) \; \vert \, x \in M , V^{m-n} \subset T_x M \} .$$
Suppose that $\sigma (M)^*_{m-n} > 0$. The following theorem is related to integral geometry on Riemannian manifolds.

{\bf Theorem 2.1}. {\it Let $W$ be a compact $(m-n)$-dimensional submanifold in $M$. Then its volume can be
estimated from below: }
$$ \text{  vol}(W) \geq  (\sigma (M)^*_{m-n}) ^{-1} \int_{(M)^*} \#(W \bigcap N_y) \, \mu_y
                                    . \eqno(2.2)   $$

{\bf Proof.} It is easy to find a finite triangulation $W_i^{\varepsilon}$ of $W$ by simplices of diameter less than $\varepsilon$, that is, $ W
= \bigcup_i W_i^{\varepsilon}$ and $ \text{  vol}_{m-n}(W_i^{\varepsilon} \bigcap W_j ^{\varepsilon} )\: = \: 0 $  if $ i \not= j $, such that for every $i$ the number of connected components of the
intersection of $W_i^{\varepsilon}$ with any
submanifold $N_y$ is at most one. So we have:

$$ \text{  vol} (W) = \sum _i \text{  vol}(W_i^{\varepsilon}  ) , \eqno(2.3)$$

$$ \int_{(M) ^*} \# (W \bigcap N_y  )  \, dy  = \sum_i \int_{(M)^* }
\# (W_i^{\varepsilon} \bigcap N_y ) \, dy    . \eqno(2.4) $$

With the help of (2.3) and (2.4) Theorem 2.1 can be proved if we show (2.2) for $W_i^{\varepsilon}$
instead of $W$. Hence, in view of our assumption it suffices to prove:
$$ \text{  vol}(W^{\varepsilon}  ) \ge (\sigma (M)^*_{m-n} )^{-1} \int_{S_{W^{\varepsilon}} }
               \, \mu_y  . \eqno(2.2.\varepsilon ) $$

Letting $ \varepsilon \to 0 $  we get the infinitesmal version of $ (2.2.\varepsilon)$:
$$ \lim_{\varepsilon \to 0} { \text{  vol}\,  (W^{\varepsilon} ) \over \text{  vol} ( S _{W^{\varepsilon} }) }= cd(x,T_xW^{\varepsilon})^{-1} \geq (\sigma (M)^* _{m-n})^{-1}    .    \eqno(2.2.0 )  $$

Obviously, $(2.2.0 )$ follows from (2.1). By integrating we obtain
$(2.2.\varepsilon)$. This completes the proof.

In the example 2.1,
if we exhaust $TM^m$ by compact bundles $TM_R$ of tangent vectors of length
$R$ over $M$, then we can also get the deformation coefficient $\sigma (TM_R)^*_m = 1$.
Here the set $(TM_R)^* $ consisting of  $m$-dimensional tangent balls of radius $R$ is diffeomorphic to $ M$. 

{\bf Corollary 2.2}. {\it Lower bound of the volume of nontrivial cycles in
Riemannian manifolds. Suppose $N \subset M$ is a k-dimensional submanifold realizing a
nontrivial cycle $[N] \in H_k (M^{n+k}, G) ,\; G= {\mathbf Z}$ or ${\mathbf Z_2}$. Let 
$(M)^*$  be a family of  n-dimensional submanifolds $N^* _{\lambda}$
realizing a nontrivial cycle $[N^*] \in H_n (M^{n+k}, G)$.
Let $\chi$ be the (algebraic)
intersection number of $[N]$ and $[N^*]$. Then we get:                }
$$ \text{  vol} (N) \ge \chi \cdot (\sigma (M)^*_k) ^{-1} \cdot \text{  vol}(M)^*.  $$

We note that Theorem 2.1 is still valid for a compact $k$-dimensional set $W$
almost everywhere smooth except singularities of codimension 1. On the other
hand, it is well-known that homological volume-minimizing cycles
are such sets [Fe 1].  So Corollary 2.2 yields the following criterion for
global minimality.

{\bf Corollary 2.3.} {\it Let $N \subset M$ be a k-cycle almost everywhere smooth except singularities of codimension 1. Suppose that the inequality in Corollary 2.2 is   an equality for $N$. Then $N$ is a globally minimal cycle.}

{\it Example 2.2}. Consider the group $U_n$ equipped with the standard
bi-invariant
metric, that is, on the tangent space $T_eU_n = \frak {u}_n$ this metric is defined as follows: $<\xi , \eta > = -tr(\xi \eta)$. Applying Corollary 2.3 we will show that the subgroup $S^1$ of all diagonal scalar elements is a homological minimal submanifold. Indeed, $U_n$
is a fibred space over $S^1 : g \mapsto \det (g)$, whose fibers are congruent with the subgroup $SU_n$. First, we note that $SU_n$ meets $S^1$ at
exactly $n$ points $x_k = diag (\exp { 2ki\pi \over n});\; k=0,\cdots ,n-1$. Therefore,
any  fibre $a \cdot SU_n$, $a \in S^1$,  meets $S^1$ exactly n points $ a\cdot x_k$. Clearly, at every intersection point $y=a\cdot x_k$  the
tangent spaces $T_ yS^1$ and $T_y (a\cdot SU_n)$ are perpendicular. Further, we observe that the algebraic intersection number between $S^1$ and $SU_n$ 
equals $n$ since $S^1$ is homologous to n times of the circle $U_1$ which
generates the homology group $H_1(U_n,{\mathbf Z})$.
 Now, it is easy
to see that if we set $(M)^*$ to be the space of cosets of the subgroup $SU_n$
in $U_n$, then $\sigma (M)^* _1 =1$, and by Corollary 2.3, $S^1$
has the minimal length in its homology class of $H_1 (U_n,{\mathbf Z}))$.

In most of our applications we are interested in cycles of compact  homogeneous
Riemannian spaces. We shall denote $(\cdot )$ the group multiplication or the 
action of a group on homogeneous spaces. Sometimes we omit this notation
 $(\cdot )$ if no confusion arises. Let $M =G/H$, where   $H$ is a compact 
subgroup in a compact group
$G$. Let $K$ be another compact subgroup of $G$. Denote $L$ the intersection of
$H$ and $K$. We consider the space $(M)^*$ of all submanifolds $g \cdot K/L \subset  G/H$
which are obtained from $K/L$ by the left shift $g, g\in G$. Obviously, $G$ acts transitively
on $(M)^*$. Let us denote $I(K)$ its isotropy group at the point $e \cdot K/L \in (M)^*$.

{\bf Lemma 2.4}. {\it The isotropy group $I(K)$ coincides with the subgroup
$K \cdot (H \cap N(K))$, where   $N(K)$ is the normalizer of the subgroup $K$ in $G$.}

{\bf Proof}. Clearly, the subgroup $I(K)$ consists of all elements $g \in G$ such that
$g \cdot K \subset K \cdot H$. So we have
$$ I(K) = \bigcap_{k \in K} (K \cdot H \cdot k) = \bigcup_{h \in H} \{ \bigcap_{k
\in K} (K \cdot h \cdot k) \} . $$
Let $h \in H$ be an element such that the intersection $\bigcap_{k \in K} (K \cdot
h \cdot k)$ is not empty. We easily verify that the last condition is equivalent
to $h$ being an element of the
normalizer $N(K)$. Hence the lemma follows immediately.

The condition under which submanifold $y \cdot K/ L \subset M$ contains a point 
$x = (g\cdot H) / H \in M$ is the relation $ y \in g\cdot H \cdot K$. So
we have the following lemma.

{\bf Lemma 2.5}. {\it Let $x = \{ gH \} \in M=G/H $. Then the set $S_x 
\subset (M)^*
= G/I(K)$ is the submanifold $gH/L' $, where   $L' = H \cap I(K)$.}

Our purpose now is to compute the
deformation coefficient $cd(x,V)$ for $x \in M$. Without loss
of generality we can assume that $x = \{ eH \} $, and then $V \subset T_{\{eH\}}M$. Denote by $\frak{g}$ the Lie algebra of $G$. Let us consider the map $\pi _*:  \frak{g} \to T_{\{eH\}}M$
which is induced by the natural projection $\pi: G \longrightarrow G/H =M$.
Let $\frak{h^g}$ be the orthogonal complement (with respect to some $Ad_G$-invariant metric on $\frak{g}$) to the subalgebra $\frak{h}$ in
$\frak{g}$. Then we identify $T_{\{eH\}}M$ with $\frak{h^g}$ by the
map $\pi _*$. This isomorphism $\pi _*$ is an isomorphism of $Ad_H$-modules. From now on we consider the metric on $\frak{h^g}$ which is induced by the isomorphism $\pi _* ^{-1}$.

{\bf Proposition 2.6}. {\it Let} $k= \text{codim} (K/L)$. {\it Then the k-dimensional deformation
coefficient $ cd( \{ eH\} , V^k)$ depends only on the H-action orbit passing through the k-dimensional
subspace $V^k$ on the space $\bigwedge ^k (\frak {h^g})$ .}

{\bf Proof}. Let us denote by $\exp$ the exponential map from Lie algebra onto Lie group. We note that we can replace the family of exhausting geodesic balls $B(\{ eH\} ,V,r)$ and
the corresponding set $S_{B(\{ eH\} ,V,r)}$ in the formula (2.1) by any family
of exhausting submanifolds $B'(\{ eH\} ,V,r)$ and $S_{B'(\{ eH\} ,V,r)}$ such that: $T_{\{ eH\} } 
B'(\{ eH\} ,V,r) = V$, $B'(\{ eH\} ,V,r_1) \subset B'(\{ eH\} ,V,r_2)$ if $r_1 \le r_2$, and
$B'(\{ eH\} ,V,r) \to \{ eH\}$ when $r \to 0$. We choose $B'(\{ eH\} , V,r)= \{\exp V(r) \cdot H\} /H$, where $V(r)$ denotes the
ball of radius $r$ in the tangent space $V \subset \frak{h^g} \subset \frak{g}$. Hence, according to Lemma 2.5 we get $S_{B'(\{ eH\} ,V,r)} = \exp V(r) \cdot H/L'$. Therefore we obtain
$$ cd( \{ eH \} , V) = \lim_{ r \to 0} { \text{vol} (\exp V(r) \cdot H/L') \over
                       \text{vol}( \exp V(r)\cdot e/H ) }  .\eqno(2.5)$$

We choose an orthonormal basis of vectors $\{ v_i \}$ in $V$.
Fix a point $ x =\{ \tilde{x}L' \} \in H/L' \subset G/I(K)$, where   $\tilde{x} \in H \subset G$. The tangent space to
$\exp V(r) \cdot H/L'$ at the point $x$ is the sum of the tangent spaces $T_x(H/L')$ and
$T_x(\exp V(r)\cdot x)$. Consider the map
$$\rho :V(r) \longrightarrow \exp V(r) \cdot \{\tilde{x}L'\} ; \hspace{1cm} v \mapsto \exp v \cdot \{\tilde{ x}L\}' .$$

Its differential $d\rho$ sends the vector
$v_i$ to the projection of the vector 
$${d\over  dt} \exp t  v_i \cdot \tilde{x}|_{t=0} \in T_x G$$
on the tangent space  $T_x (G/I(K))$ since $G/I(K)$ is
the quotient space of the right $I(K)$-action on $G$ . Denote $\hat{v_i} (x)$ the
resulting vector $d\rho (v_i) \in T_x (G/I(K))$. Then we have $T_x(\exp V(r)\cdot x ) =\text{  span}\{ \hat{v}_i, i= 1, \cdots ,n \}$. So (2.5) can be rewritten as follows:
$$ cd( \{ eH \}, V) =  \int_{H/L'} \text{  vol} ( \overline{T_x (H/L')} 
\wedge \hat{V}_x ) \,
                       \mu_x   ,\eqno(2.6)  $$ 

where $\overline{T_x (H/L')}$ denotes the unit polyvector associated with $T_x (H/L')$, and $\hat{V}_x = \hat{v_1}(x) \wedge ... \wedge \hat{v_k}(x)$. First, we note that $ \text{  vol} ( \overline{T_x (H/L')} 
\wedge \hat{V}_x ) = | < \hat{V}_x, \overline{W_x}>|$, where  the associated subspace $W_x$ is the orthogonal complement to $T_x (H/L')$ in  $T_x(G/I(K))$.  Secondly, we observe that for each $h \in H$ we have
$$ \widehat{Ad_h v_i} (x) =  h_* \hat{v}_i(h^{-1} \cdot x).$$
Therefore we obtain
$$ cd( \{ eH \}, Ad_h V) = \int_{H/L '} |< h \cdot \hat{V}_{h^{-1}\cdot x}, \overline{W_x}>|\, \mu_x .  \eqno(2.7)$$
Now Proposition 2.6 immediately follows from (2.6), (2.7) and the $G$-invariance of the metric on $G/I(K)$.

Let us consider the case when the invariant metrics on $G/H$ and $G/I(K)$ are canonical (i.e. thay are obtained from a bi-invariant metric on $G$ factorized by the action of its subgroups $H$ and $I(K)$ respectively.) In this case the formula for $cd(\{ eH\},V)$ has a very simple expression. Denote by $\frak{h}$ and $\frak{k}$ the Lie algebras of the subgroups $H$ and $K$ respectively. Let  $W$ be the orthogonal complement to the span of these subalgebras in $\frak{g}$, that is,
$$ \frak{g} = W \oplus ( \frak{h} + \frak{k}).$$
 Then we obtain the following lemma.

{\bf Lemma 2.7}. {\it Under the above assupmtions we have}
$$ cd( \{ eH \}, V) = \int_{H/L '} | <\overline{V}, Ad_{\tilde{x}} (\overline{W})> | \, dx .\eqno(2.8)$$
{\it Proof}. Denote by $pr( \frak{h})$ the orthogonal projection of $\frak{h}$ onto the orthogonal complement to $\frak{k}$ in $( \frak{h} + \frak{k})$. We have the following orthogonal decomposition
$$T_{\{ eI(K)\} } G/I(K) = W \oplus \text{  span}\{ z \in pr(\frak{h}) | 
\; <z, \frak{h}\cap \frak{n(k)}> = 0 \; \} ,$$
and
$$T_{\{ eL'\} } H/L' = \text{  span}\{ z \in pr(\frak{h}) | \; <z, \frak{h}\cap \frak{n(k)} = 0\; \} .$$
Therefore, the normal fiber $W_{\{ eL'\}}$ coincides with $W$. Since $\tilde{x} \in H$ the shift $L_{\tilde{x}}$ preserves the normal bundle of $H/L'$ in $G/I(K)$. Hence, $W_x$ = $\tilde{x}_* W$. 

Our next
aim is to compute $\hat{v_i}(x)$. Let us choose an orthonormal basis $ f_1..., f_N $
of the space $lI(K)^G = T_{\{ eI(K) \} } G/I(K) $. The shift $L_{\tilde{x}} : G/I(K) \longrightarrow G/I(K) , \; \{ gI(K) \}
\mapsto \tilde{x} \cdot \{ gI(K) \} $, sends the vector $f_i$ to the vector $f_i ^{\tilde{x}}(x)$ . Obviously,
$f_i ^{\tilde{x}}(x)$ is an orthogonal basis of the tangent space $T_x( G/I(K))$. Straightforward calculation shows that 
$$ < \hat{v_i}(x) , f_j ^{\tilde{x}}(x) > \; = \; <v_i , Ad_{\tilde{x}} f_j  >, $$
where $<,>$ in the right hand side of the above formula denotes the restriction of the bi-invariant metric on $G$ to the algebra ${\frak g}$.

Now, taking into account  (2.7) (with $h = e$) we immediately get the formula (2.8). Clearly, the
space $W$ is invariant under the action $Ad_{L'}$. Thererfore,
the integrand on the right hand side of (2.8) depends only on $x$. This
completes the proof of Lemma 2.7.

{\it Example 2.3}. Let $M =S^n =SO_{n+1}/SO_n$, and $(M)^* =SO_{n+1}/S(O_{k+1}\times O_{n-k})$
the set of great (totally geodesic)
$k$-dimensional spheres in $S^n$. Here $H = SO_n$ acts on the Grassmannian $G_{n-k}
(T_e M) \cong SO_n /S(O_k \times O_{n-k})$ transitively. This means that $cd(x,V)$ is a constant $\zeta _{n-k}$. Taking into account ($2.2.\varepsilon$), (2.2.0)
(which become  equalities in this case) and (2.3), (2.4) we get:

{\bf Proposition 2.8} [Sa]. {\it Let $N^{n-k}$ be a submanifold in $S^n$. Then its volume
can be computed from the following formula:}
$$ \text{vol}(N^{n-k}) = \zeta_{n-k} \cdot \int_{SO_{n+1}/ S(O_{k+1} \times  O_{n-k})} \# (N^{n-k} \bigcap S^k (x))
\, \mu_x , $$

{\it where } $\zeta_{n-k} = 1/2 \, \text{vol}(S^{n-k}) \cdot \text{vol}(SO_{n+1}/S(O_{k+1} \times
 O_{n-k} ))^{-1} $.

The same formula holds for a submanifold $N^{n-k} \subset {\mathbf R}P^n$, but we should replace
$S^k$ by ${\mathbf R}P^k$. Further, we note that any projective space ${\mathbf R}P^k$ meets almost
all projective spaces of complementary dimension at one point (cf. Proposition
3.6). Hence in view of Corollary 2.3 we obtain:

{\bf Proposition 2.9}. {\it The projective space ${\mathbf R}P^k$ has the minimal volume
in its homology class $[{\mathbf R}P^k] \in H_k ({\mathbf R}P^n ,{\mathbf Z_2}) ={\mathbf Z_2} $. }

This proposition was obtained by Fomenko [Fo 1] using a different method of geodesic defects.
 
{\it Example 2.4}. Let $M = {\mathbf C}P^n = U_{n+1} /(U_n \times U_1) $. Then $T_e {\mathbf C}P^n =
{\mathbf C}^n= {\mathbf R}^{2n}$, and $H = U_n \times U_1$ does not act on $G_k({\mathbf R}^{2n})$ transitively.
But $H$ acts on the complex Grassmannian $G_k ({\mathbf C}^n)$ transitively, and $H$ also acts
on the Lagrangian Grassmannian $GL({\mathbf C}^n) = U_n/O_n$ transitively. Considering the
family $(M)^* _1 = U_{n+1} / (U_{n-k+1} \times U_k ) $ of all canonically embedded complex projective spaces of
dimension $(n-k)$  in $M$, and the family $(M)^* _2 = U_{n+1}/O_{n+1}$ of all canonically 
embedded real projective spaces of dimension $n$ in $M$, we get:

{\bf Proposition 2.10}. {\it a) Crofton type formula. Let $N^{2k}$ be a complex manifold
in ${\mathbf C}P^n$. Then its volume can be computed  from the following formula:}
$$ \text{vol}(N^{2k}) = \zeta ^C _k \cdot \int_{U_{n+1} /( U_{n-k+1} \times U_k ) } \# (N^{2k} \bigcap {\mathbf C}P^{n-k}(x)) \,
                 \mu_x ,    $$

{\it where the constant $\zeta ^C _k$ does not depend on $N^{2k}$.

b) Let $N^n$ be a Lagrangian manifold in ${\mathbf C}P^n$. Then its volume can be computed
from the following formula:}
$$ \text{vol}(N^n) = \zeta ^L _n \cdot \int_{U_{n+1}/O_{n+1}} \# ( N^n \bigcap {\mathbf R}P^n (x)) \,
               \mu_x ,  $$ 

{\it where the constant $\zeta ^L _n$ does not depend on $N^n$ and $U_{n+1} / O_{n+1}$
is the space of all real projective spaces of dimension $n$ in ${\mathbf C}P^n$. }

When $k =1$ we have the following inequality.

{\bf Proposition 2.11}. {\it Integral Wirtinger Inequality. Let $N^2$ be a real surface in ${\mathbf C}P^n$. Then the following inequality holds
$$ \int_{{\mathbf C}P^n} \# (N^2 \bigcap \gamma) \, d \gamma  \le  Area (N^2) , $$

where $\gamma$ is a (complex) hyperplane of ${\mathbf C}P^n$, and the space of these hyperplanes
is identified with ${\mathbf C}P^n$ equipped with the invariant measure. Moreover,
the inequality becomes   an equality if and only if $N^2$ is a complex
curve.}

{\bf Proof}. We consider the family $({\mathbf C}P^n)^*$ of complex hyperplanes
in ${\mathbf C}P^n$. According   to Theorem 2.1 it suffices  to show  that the associated deformation coefficient  $cd(x,V^2)$ attains its  maximal value if   and  only if $V^2$ is  a complex line. Using the above notations  we have $H= U_n \times U_1$, $K= U_1 \times U_n$, $L= L'= U_1 
\times U_{n-1} \times U_1 $, and then $H/L' = {\mathbf C}P^{n-1}$. With the help of (2.8) we get
$$cd(\{ eH\} , V^2) =  \int_{{\mathbf C}P^{n-1}}   | <\overline{V^2}, Ad_{\tilde{x}} (\overline{W})>|\, dx .$$
Let $L'' = \{ 1 \} \times U_{n-1} \times U_1$. Then $S^{2n-1} = H/L''$
is also considered as the unit sphere in the orthogonal complement
$(\frak{l}'')^H$ to $\frak{l}''$ in $\frak{h}$.
We consider the Hopf fibration $S^{2n-1} \longrightarrow {\mathbf C}P^{n-1}$. It is well-known
that the Hopf fibres are the $U_1$ orbits on $S^{2n-1}$, and the invariant Riemannian  metric on ${\mathbf C}P^n$ is obtained from the one on $S^{2n-1}$ factorized  by the $U_1$ action. Therefore we get
$$ cd(\{ eH\} ,V^2) = \text{  vol} (U_1)^{-1} \int_{S^{2n-1}} |<\overline{V^2}, Ad_{\tilde{x}} (\overline{W})>|\, dx .$$

Now we apply the normal form theorem of Harvey and Lawson to $V^2$.

{\bf Lemma 2.12}. [H-L, Lemma 6.13]. {\it There exists a unitary basis $v_i, Jv_i$ in ${\mathbf C}^n = T_{\{ eH\} } {\mathbf C}P^n$  such  that $V^2= \cos \tau \cdot v_1 \wedge Jv_1 + \sin \tau \cdot v_1 \wedge v_2$}.

Taking into account the equality $Ad_{\tilde{x}} (\overline{W}) = x \wedge Jx$ for $x \in S^{2n-1} \subset (\frak{(l'')^h}$ we obtain
$$ cd(\{ eH\} ,V^2) = \text{  vol} (U_1)^{-1}\int_{S^{2n-1}}|<\cos \tau \cdot v_1 \wedge Jv_1 +
\sin \tau \cdot 
 v_1 \wedge v_2  , x\wedge Jx>|\, dx .\eqno(2.9)$$
Let $a_i(x)= <x,v_i>$ and $b_i(x)= <x,Jv_i>$. From (2.9) we get
$$cd(\{ eH\} ,V^2)= \text{  vol}(U_1)^{-1} \int_{S^{2n-1}} | (a_1 ^2 (x)+ b_1 ^2(x))\cos \tau + 
(-a_1(x)b_2(x)
+ $$
$$+a_2(x) b_1(x))\sin \tau |\, dx .\eqno(2.10)$$
Since the integrand in (2.10) is homogeneous of degree 2 on
${\mathbf R}^{2n}$, we observe that our calculation
can be reduced to the one on sphere $S^3$. Namely, there exists a constant $\chi _n$ such that 
$$cd(\{ eH\} ,V^2) = \chi _n \int_{S^3} |(a_1 ^2 (x)+ b_1 ^2(x))\cos \tau + \sin \tau
(-a_1(x)b_2(x)+ a_2(x) b_1(x))\sin \tau |\, dx . $$

Hence we obtain
$$cd(\{ eH\} ,V^2) \le \chi _n (\int_{S^3}  | a_1 ^2 (x)\cos \tau -  a_1 (x) b_2 (x)\sin \tau |\,$$
$$ + | b_1 ^2 (x)\cos \tau  + a_2 (x)b_1 (x)\sin \tau |  \, dx )\eqno(2.11)$$
We choose the torus coordinates on $S^3$. Namely we put 
$$a_1 (x) = \sin \beta (x)\cos \alpha (x),\; a_2 (x) =\sin \beta (x)\sin \alpha (x),$$
$$b_1 (x) = \cos \beta (x) \cos \gamma (x),\; b_2 (x) = \cos \beta (x)\sin \gamma (x),$$
where $\beta \in [0,\pi], \alpha \in [0,2\pi], \gamma \in [0,2\pi]$. So, the
action of the group $S^1 \times S^1$ on $S^3$ given by:
$\alpha (x) \to \alpha (x) + \theta _1, \; \gamma (x) \to \gamma (x) + \theta _2$
preserves the invariant measure on $S^3$. In these coordinates (2.11) becomes  
the following inequality
$$cd(\{ eH\} ,V^2)\le \chi _n (\int_{S^3} |\sin ^2 \beta \, \cos  \alpha \, \cos (\alpha + \tau)|\, \mu (\alpha , \beta , \gamma )+$$
$$+ \int_{S^3} |\cos ^2 \beta \, \cos \gamma \, \cos (\gamma - \tau)|\, \mu 
(\alpha , \beta , \gamma ) ,\eqno(2.11')$$
where $\mu $ is the invariant measure on $S^3$.
Applying the Schwarz  inequality for integrals to the right hand side of $(2.11')$ we get
$$cd(\{ eH\} ,V^2) \le \chi _n \{ (\int_{S^3} |\sin ^2 \beta \, \cos ^2 \alpha |\, \mu )^{1/2} \cdot
(\int_{S^3} |\sin ^2 \beta \, \cos ^2 (\alpha +\tau )|\, \mu )^{1/2} + $$
$$+ (\int_{S^3} |\cos ^2
\beta \, \cos ^2 \gamma |\, \mu )^{1/2} \cdot (\int_{S^3} |\cos ^2 \beta \, \cos ^2 (\gamma - 
\tau )|\, \mu )^{1/2}) . $$
As it was mentioned above the transformation  $g(\tau )$: $\alpha \to \alpha +\tau, \; 
\gamma \to \gamma -\tau$ , preserves the invariant measure $\mu$. Therefore we get
$$cd(\{ eH\} ,V^2) \le \chi _n \int_{S^3} |\sin ^2 \beta \cos ^2 \alpha + \cos ^2 \beta \cos ^2
\gamma |\, \mu . \eqno(2.12)$$
The inequality  (2.11) becomes   an equality if and only if $\tau = 0$. Observe
that the right hand side of (2.12) equals $cd(\{ eH\} ,v_1 \wedge Jv_1)$.
This means that the deformation coefficient $cd(\{ eH\}, V^2)$ attains its maximal value  only at complex lines. Our proof is completed.

{\it Remark}. From the above proof we immediately deduce  a dual proposition which
replaces a two-dimensional surface $N^2 \subset {\mathbf C}P^n$ by a surface of codimension 2. A proof for the case of an arbitrary $k$ will be given in \S 3 (see Proposition 3.10 ).

\vspace{0,5cm}
 
\S 3. {\bf Minimal cycles in Grassmannian manifolds. }

We denote $G_k({\mathbf R}^n)$ the Grassmannian of unoriented
k-planes through the origin in ${\mathbf R}^n$ and its 2-sheeted covering
by $G^+ _k ({\mathbf R}^n)$. We denote $G_k({\mathbf C}^n)$ and $G_k ({\mathbf H}^n)$ the complex Grassmannian
and the quaternionic Grassmannian respectively.
The question of finding and classifying  globally minimal
cycles in Grassmannian manifolds has attracted attention of many
mathematicians. The first non-trivial result was obtained by A. T. Fomenko
in 1972 using his method of geodesic defects [Fo 1, Le-Fo] and by M. Berger
in the same year using calibration method [Be]. In particular, Fomenko proved 
 that the canonically embedded real projective space ${\mathbf R}P^l \longrightarrow
{\mathbf R}P^n, \, l \le n $, is globally minimal, and Berger proved that
${\mathbf H}P^k$ is homologically volume-minimizing in ${\mathbf H}P^n$ if $k \le n$. Recently, employing Euler forms and their "adjusted powers"
as calibration H. Gluck, F. Morgan and W. Ziller proved that {\it if k = even $\ge 4$, then each
$$G^+ _1 ({\mathbf R}^{k+1}) \subset G^+ _2({\mathbf R}^{k+2}) \subset \cdots \subset G^+ _l ({\mathbf R}^{k+l})$$
is uniquely volume minimizing in  its homology class} [G-M-Z]. H. Tasaki showed that the same
proof implies that $G_k ({\mathbf H}^{m+k})$ is uniquely volume minimizing in its homology class in
$G_n ({\mathbf H}^{m+n})$ for all $m$ , even and odd [T]. In this section using our method we
prove:

{\bf Theorem 3.1}. {\it The canonically embedded  real Grassmannian submanifold
$G_k({\mathbf R}^{k+m})$ in $ G_l({\mathbf R}^{l+m}), \; k \le l,$ has the minimal volume in
its homology class with coefficients in ${\mathbf Z} $ or ${\mathbf Z_2}$.}

We will show in \S 4 that this theorem implies the G-M-Z Theorem mentioned above.
But the G-M-Z Theorem implies our Theorem only in the case when $m$ is even and $G=
{\mathbf Z}$, because
when $m$ is odd, each $G^+ _k ({\mathbf R}^{k+m})$ bounds over the reals in $G^+ _l ({\mathbf R}^{l+m})$.

{\bf Theorem 3.1'.} {\it Classification Theorem. Let $M$ be a volume-minimizing
cycle of the non-trivial homology class $[G_k({\mathbf R}^{m+k})] \in H_* (G_l ({\mathbf R}^{m+l}), G)$,
where $G ={\mathbf Z}$ or ${\mathbf Z_2}$. Then $M$ must be one of these sub-Grassmannians.}

{\bf Theorem 3.2.} {\it The canonically embedded complex Grassmannian submanifold
$G_k ({\mathbf C}^{k+m})$ in $G_l({\mathbf C}^{l+m}), \; k \le l,$ has the minimal volume in its
homology class with coefficients in ${\mathbf Z_2}$.}

{\bf Theorem 3.3.} {\it The canonically embedded quaternionic Grassmannian 
submanifold $G_k({\mathbf H}^{k+m})$ in $G_l ({\mathbf H}^{l+m}), \; k \le l,$ has the minimal
volume in its homology class with coefficients in ${\mathbf Z_2}$.}

{\it Remark}. Of course, we can also prove these theorems with respect to integral homologies (and then real homologies) by the same method.

{\bf Proof of Theorem 3.1}. We apply results of Section 2 to $G= SO_{l+m}$, $H= S(O_l \times O_m)$,  $I(K) = K = S(O_k \times O_{l-k+m})$,
$L=L' = S(O_k \times O_{l-k} \times O_m)$. We consider the family $(M)^*= SO_{l+m}/S( O_k\times O_{l-k+m})$ 
of homogeneous subspaces obtained
from $G_{l-k}({\mathbf R}^{l-k+m})$ by the action of the group $SO( {\mathbf R}^{l+m})$ (see \S 2). Let $V$ be a $km$-dimensional subspace of $T_e G_l
({\mathbf R}^{l+m})$, where $e = \{eH\}$.
According to Lemma 2.7 we get:
$$ cd( e, V) = \int_{S(O_l  \times O_m)/ S(O_k \times
O_{l-k} \times O_m )}
       \vert < \overline{V}, Ad_{\tilde{x}}\overline{ W} > \vert \, dx  $$

$$ = \int_{SO_l /S(O_k \times O_{l-k})} | < \overline{V}, Ad_{\tilde{x}} \overline{W} > | \, dx  , \eqno(3.1) $$

where $W$ denotes the tangent space  $T_eG_k ({\mathbf R}^{k+m})$.

Clearly, the group $SO_l$ acts on the tangent space $T_e G_l ({\mathbf R}^{l+m}) =  {\mathbf R}^l \otimes
{\mathbf R}^m $ as the sum of $m$ irreducible representations $\pi _1$ of dimension $l$.
Namely, in the matrix representation of $T_e G_l ({\mathbf R}^{l+m}) \longrightarrow
so_{l+m}$ these irreducible spaces can be chosen as $m$ columns ${\mathbf R}^l _i$. Let us denote
by $I$ the canonical operator of the  decomposition $T_e G_l ({\mathbf R}^{l+m}) = \bigoplus {\mathbf R}^l _i$ with respect to the adjoint action of $SO_l$, that is, $I\cdot Ad =
Ad \cdot I$ and $I({\mathbf R}^l _i) = {\mathbf R}^l _{i+1}$. Obviously,  we have $W = W_1 \oplus
I(W_1)\oplus \cdots \oplus I^{m-1}(W_1)$, where   $W_1 = W \cap {\mathbf R}^l _1$. So we get
$Ad_g W = Ad_g W_1 \oplus I(Ad_g W_1)\oplus \cdots \oplus I^{m-1} (Ad_g W_1)$.
Now we consider the following fibration $j: SO_l /SO_{l-k} \longrightarrow SO_l /S(O_k \times 
O_{l-k})$, where   the total space is considered as the Stiefel manifold of frames of $k$ orthonormal vectors in ${\mathbf R}^l _1$, and the base is the Grassmannian of unit simple $k$-vectors in ${\mathbf R}^l$, which is identified with the set of all $Ad_{\tilde{x}} W$. Thus, if $x$ is a frame of $k$ orthonormal vectors $(v_1, \cdots , v_k)$, then $j(x) =
v_1 \wedge \cdots \wedge v_k$. Let the metrics on the above spaces be the standard ones. Since the volume of each fibre $O_k$ is a constant $\lambda _{k,l}$, we can rewrite
integral (3.1) as follows
$$cd(e ,V) = \lambda _{k,l} \int_{SO_l /SO_{l-k} } |<\overline{V}, j(x)\wedge
I(j(x))\wedge...\wedge I^{m-1}(j(x))>| \, dx .\eqno(3.2)$$
 We consider the fibration $SO_l/SO_{l-k} \rightarrow SO_l/SO_{l-k+1}$ with fibre $S^{l-k}$;
 it maps a $k$-frame $x=(v_1,\cdots,v_k)$ to a $(k-1)$-frame $x'=(v_1,\cdots,v_{k-1})$. Denote ${\mathbf R}^{l-k+1}(x')$ the linear subspace associated with the fiber $S^{l-k}$ over the point $x'$. Using integration along fibres we deduce from (3.2)
$$cd(e ,V) = \lambda _{k,l} \int_{SO_l/SO_{l-k+1}} \int_{S^{l-k} (x')}
|<\overline{V}, j(x',y)\wedge...\wedge I^{m-1}(j(x',y))>| \,dy \,dx' $$
$$= \lambda _{k,l} \int_{SO_l / SO_{l-k+1}} \{ |<V, j(x')\wedge ...\wedge I^{m-1}(j(x')) >| \cdot $$
$$\cdot \int_{S^{l-k}(x)} |< V^{\perp} (x'), y\wedge ...\wedge I^{m-1} (y)>|\, dy \} \, dx',
\eqno(3.3)$$

where $|<V, z>|$ denotes the volume of the orthogonal projection of a simple  polyvector
$z$ on the plane $V$; and $V^{\perp}(x')$ is the intersection of $V$ with the space
${\mathbf R}^{l-k+1} (x') \oplus  \cdots \oplus  I^{m-1} ({\mathbf R}^{l-k+1} (x'))$.

{\bf Proposition 3.4}. {\it Let $p \leq q$. For each $mp$-plane $V \subset {\mathbf R}^q \oplus  \cdots \oplus  I^{m-1}
({\mathbf R}^q)$, where ${\mathbf R}^q \subset {\mathbf R}^l _1$, we put
$$M(V) = \int_{S^{q-1}} |<V, x\wedge \cdots \wedge I^{m-1} (x)>|\, dx .$$
Then $M(V)$ reaches its maximal value if and only if $V = V^p \wedge \cdots  \wedge I^{m-1}(V^p)$, where   $V^p \subset {\mathbf R}^q$. }

Repeating the reduction process (3.3) and applying Proposition 3.4 we obtain
the following proposition immediately.

{\bf Proposition 3.5}. {\it The deformation coefficient cd(e,V) attains  its maximum
at $V_0$ if and only if there exists $\tilde{x} \in SO_l$ such that $V_0 = Ad_{\tilde{x}} W$.}

{\bf Proof of Proposition 3.4}. Obviously, we have
$$M(V) \le \int_{S^{q-1}} |<V, x>| \cdot ... \cdot |<V, I^{m-1}(x)>| \, dx. \eqno(3.4)$$
Applying the theorem about geometric and arithmetic means we infer from (3.4)
$$ M(V) \le ({1\over m})^{m/2} \int_{S^{q-1}} (\sum_r |<V, I^r (x)>|^2 )^{m/2} \, dx .\eqno(3.5)$$
Now we study the projection $I_V ^r (x)$ of $I^r (x)$ on $V$ and its length  $|<V, I^r (x)>|$. 
Let $B_r$ denote the symmetric bilinear form on ${\mathbf R}^q$ defined by $B_r (x,x) =< I^r_V (x), I^r _V(x)>$.
Let $\theta ^r _j$ be the eigenvalues of $B_r$, $j=1, \cdots ,q$. Evidently, $0 \le \theta ^r _j \le
1$.

{\bf Lemma 3.6}. {\it The following identity holds
$$ \sum_{r,j}\theta ^r _j = \sum_r tr(B_r) =\dim V = mp .$$}
{\bf Proof}. Let   $\Pi _r$ be the bilinear
form on $V$ defined by
: $\Pi _r (x,x) = < \pi _r (x), \pi _r (x)>$, where   $\pi _r$
denotes the orthogonal projection on $I^r ({\mathbf R}^q)$. We will show that $tr(B_r) = tr(\Pi _r)$. Without loss of generality we can assume that
$\dim V \ge \dim I^r({\mathbf R}^q)$. Now we consider
the eigenvectors $\{ f_i^r \} \in I^r ({\mathbf R}^q)$ of $B_r$ corresponding to $\theta ^r _i$.
Then $\{ f_i^r\}$ can be chosen as an orthonormal basis in  $I^r({\mathbf R}^q)$. Clearly, we have
$$ <f_i^r, I^r _V(f_j^r)> = <I^r _V (f_i^r), I^r _V (f_j^r)> = \delta _{ij} \theta ^r _i .\eqno(3.6) $$
 We want to find the orthogonal projection $\widehat {I^r_V(f^r_j)}$
of the vector $I^r_V (f^r_j)$ $\in V$ 
on $I^r({\mathbf R}^q)$. We note that this projection is defined  uniquely,
up to multiplication by a constant, by the hyperplane orthogonal to it in the subspace $I^r({\mathbf R}^q)$. Obviously, this hyperplane $H^r _j$ is defined by the 
following equation 
$$H^r_j = \text{  span} \{ z | < I^r_V (f^r_j), z> = 0 \; \}.\eqno (3.7)$$
 Now, comparing (3.7) with (3.6), it is easy to see that $ \widehat{I^r_V( f^r_j)}\in \text{  span}\{f^r_j\}$. Therefore,
the orthogonal projection of the vector $I^r _V (f^r _j)/|I^r _V (f_j^r)| \in V$ on the subspace $I^r ({\mathbf R}^q)$ is $\theta ^r _j \, f^r_j$. Note that for any vector $w$ in the orthogonal complement to $\text{  span}
\{ I^r _V (f^r_j) \}$ in $V$ we have $<w, f^r_i> =0$. Hence, in view of (3.7), we have that $\theta ^r _j$, $j=1, \cdots ,q$, and $0$ with multiplicity $mp -q$ are eigenvalues of $\Pi _r$, and then we have $tr(B_r) =tr(
\Pi _r)$.

Further we note that $\sum \Pi _r(x,x) = <x,x>$. Therefore $\sum tr(B_r ) = \sum tr(\Pi _r) =
\dim V$. This completes the proof of Lemma 3.6. 

Let us continue the proof of Proposition 3.4. From the proof of Lemma 3.6 we know that
$$\sum_{r=0}^{m-1} |<V, I^r (x)>|^2 = \sum_{r=0} ^{m-1} B_r(x,x).$$ 
We set $B(x,x) = \sum B_r(x,x)$. Since $B_r(x,x)$ are symmetric bilinear forms whose eigenvalues belong to the
segment $[0,1]$, the symmetric bilinear form $B(x,x)$ is also positive, moreover,
its eigenvalues belong to the segment $[0,m]$. Denote these eigenvalues by
$\eta _i$, $i=1,\cdots ,q$. From Lemma 3.6 we know that $\sum \eta _i = Tr(B) = \sum Tr(B_r)=
\dim V = pm$. Let $w_i$ be the eigenvectors corresponding to $\eta _i$. 
Obviously, we can choose $w_i$ as an orthonormal basis in ${\mathbf R}^q$.
So, we rewrite (3.5) as follows

$$M( V) \le ({1\over m})^{m/2} \int_{S^{q-1}}(\sum_j \eta  _j (x _j)^2)^{m/2}\,
dx , \eqno(3.8) $$
where $x_j$ is the $j^{th}$ coordinate of $x \in S^{q-1}$ with respect to the basis
of vectors $\{ w_i \}$.
Let $F(\eta _1,\cdots ,\eta _q)$ be the function in the right hand side of (3.8) whose variables satisfy the following condition:
$$\eta _i \in [0,m]; \; \sum \eta _i = mp  . \eqno (C) $$
We want to find the maximum of $F$. To see this we choose any two variables $\eta _1$ and $\eta _2$
among $\eta _j$ and fix the others. So, we have $\eta _2 = c- \eta _1$, where     c is some constant. 
Straightforward calculation yields:
$${d^2 \over d\eta _1 ^2}( F )= ({1\over m})^{m/ 2} \int_{S^{q-1}}({m\over 2} -1) \cdot {m\over 2}
\cdot \{ \sum_j \eta _j (x _j)^2 \} ^{{m-4\over 2}} \cdot (x_1 ^2 - x_2 ^2)^2 \, dx  .$$
If $m \ge 3$ the above formula shows that $F$ is a convex function with respect to $\eta _1$. Therefore, $F$ attains its maximal value at only "boundary" variables. This means 
that under the condition C we have
$$F (\eta _1,\cdots ,\eta _q) \le F (m,...,m,0...,0) . $$
This formula shows that $M(V)$ attains its maximal value if and only if the eigenvalues of $B(x,x) = \sum B_r (x,x)$ are $(m,...,m,0,...0)$. Since $\theta ^r _j 
\in [0,1]$ we immediately obtain that for every $r$ the eigenvalues of $B_r$
are $(1,...,1,0,...0)$, moreover $B_i =B_j$ for all $i,j$. Consequently, we have
$V = V_1 \wedge I(V_1)\wedge...\wedge I^{m-1}(V_1)$. If $m =2$ then $F$ is a linear function with respect to $\eta _j$. In this case it suffices to consider two
inequalities (3.4) and (3.5) to obtain our assertion. This completes the proof
of Proposition 3.4.

Now we study the intersection between Grassmannian submanifolds in $G_l ({\mathbf R}^{l+m})$.

{\bf Proposition 3.7}. {\it For almost all (in dimension sense) $y \in (M)^* =
SO_{l+m}/S(O_k \times O_{l-k+m} )$ the space $N_y = \tilde y \cdot G_{l-k}({\mathbf R}^{l-k+m})$ meets $G_k ({\mathbf R}^{k+m})$
at only one point.}

{\bf Proof}. Geometrically, the embedding $ G_k({\mathbf R}^{k+m}) \longrightarrow G_l 
({\mathbf R}^{l+m})$ can be described as follows:
 $$G_k( {\mathbf R}^{k+m}) \ni x \mapsto x \wedge v_{l-k} \in G_l ({\mathbf R}^{l+m}) ,$$
 where $v^{l-k}$ denotes the subspace orthogonal to ${\mathbf R}^k$
in ${\mathbf R}^l$. So, the intersection $T(y)$ of the considered Grassmannians consists
of those $l$-dimensional subspaces $W^l$ such that:
$$ W^l \in ( G_k ({\mathbf R}^{k+m}) \wedge v^{l-k}) \bigcap (G_{l-k}( {\tilde y} \cdot
{\mathbf R}^{l-k+m} )\wedge {\tilde y} \cdot v^k ) .\eqno(3.9) $$

Clearly, the following lemmas yield Proposition 3.7.

{\bf Lemma 3.8}. {\it The set of all elements $y \in (M)^*$ such that
the dimension of ${\tilde y} \cdot {\mathbf R}^k \bigcap {\mathbf R}^{l-k}$ is greater than or equal to
1 has codimension 1.}

{\bf Lemma 3.9}. {\it If ${\tilde y} \cdot {\mathbf R}^k \bigcap {\mathbf R}^{l-k}$ contains
only the origin in ${\mathbf R}^{l+m}$ then $T(y)$ contains only one element.}

{\bf Proof of Lemma 3.8}. It suffices to prove that the set of
${\tilde y} \in SO_{l+m}$ with the above property has codimension 
greater than or equal to 1 in $SO_{l+m}$. Let ${\tilde y}$ belong to this set. Then
its entries (we consider ${\tilde y}$ as a matrix) satisfy the equation:

$$ \text{  vol} ( {\tilde y} \cdot v^k \wedge v^{l-k} )  = 0  . \eqno(3.10)  $$

The solution to (3.10) is an algebraic hypersurface in $SO_{l+m}$. This
completes the proof.

{\bf Proof of Lemma 3.9}. Let $W^l \in T(y)$. According to (3.9) $W^l$
contains both ${\mathbf R}^{l-k}$ and $\tilde{y} \cdot {\mathbf R}^k$. By our assumption $W^l$ must
be their span. This yields the assertion.

Let us complete the proof of Theorem 3.1. Suppose $V$ is a submanifold
of $G_l({\mathbf R}^{l+m})$ representing the same homology class as $G_k({\mathbf R}^{k+m})$. Then $V$ meets every submanifold
$N_y = {\tilde y} \cdot G_{l-k}({\mathbf R}^{l-k+m})$ at least one time. Hence, our
theorem immediately follows from Proposition 3.5, Proposition 3.7 and 
Corollary 2.3. 

{\bf Proof of Theorem 3.1'}. Let $N$ be a volume-minimizing cycle in the
homology class $[G_k({\mathbf R}^{m+k})]$. First, we observe that $N$ is almost
everywhere smooth (see [Fe 1]) and then we can apply Corollary 2.2 to $N$. On the other
hand, since $G_k ({\mathbf R}^{m+k})$ satisfies the condition in Corollary 2.3,
we conclude that the cycle $N$ also satisfies this condition. In
particular, we obtain that for almost all $x\in N$ (in dimension sense) the tangent space $T_x N$ to
$N$ satisfies the condition of maximal deformation coefficient :
$cd(x, T_x N) = \sigma (M)^* _{km}$.
In view of Proposition 3.6 we obtain that the tangent space $T_x N$ is
also tangent to some sub-Grassmannian $g\cdot G_k({\mathbf R}^{k+m})$. Then we can apply
Proposition 3.2 in [G-M-Z], which states that such a submanifold must be
one of the sub-Grassmannians $g\cdot G_k ({\mathbf R}^{k+m})$. Indeed, Proposition 3.2 in
[G-M-Z] is stated for the case of Grassmannian of oriented planes $G^+ _k ({\mathbf R}^{k+m})$,
but their Grassmannian and ours one are locally isometric, so their 
Proposition is still valid in our case. This completes the proof of
Theorem 3.1'. 

{\bf Proof of Theorem 3.2}. The proof of this theorem is similar to 
that of Theorem 3.1. First we will prove the Integral Wirtinger Inequality
for arbitrary $k$ (cf. Proposition 2.10.a).

{\bf Proposition 3.10}. {\it Let $N^{2k}$ be a manifold in ${\mathbf C}P^n$. Then its
volume can be estimated from below by}
$$ \text{vol}(N^{2k}) \ge \zeta ^C _k \cdot \int_{U_{n+1} /(U_{n-k+1} \times U_k )}
\# (N \bigcap {\mathbf C}P^{n-k} (x)) \, \mu _x ,$$
{\it where $\zeta ^C _k$ is the constant in Proposition 2.10.a. Moreover, the inequality becomes   an equality if and only if $N^{2k}$ is 
a complex submanifold.} 

{\bf Proof}. As in the proof of Proposition  2.11, it suffices to show
that the deformation coefficient $\chi _{2k} (e, V^{2k})$, related to the
family of complex projective subspaces of dimension $(n-k)$ in ${\mathbf C}P^n$,
reaches its maximal value iff $V^{2k}$ is a complex space. According to
(2.8) we obtain (see also Proposition 2.10.a):
$$ cd(e, V^{2k}) = \int_{G_k (T_e {\mathbf C}P^n) } |< \overline{V^{2k}}, Ad_{\tilde{x}} (\overline{W} )> |\, dx ,$$
where $W$ is the tangent space to the (fixed) complex projective space ${\mathbf C}P^k$.
Now we
consider the complex Grassmannian $G_{k-1} (T_e {\mathbf C}P^n)$. We associate to each point
$x \in G_{k-1} (T_e {\mathbf C}P^n)$ the fibre $q(x)$ of complex lines in the  complex $(n-k+1)$-dimensional
dimensional
orthogonal complement to the space $\text{  span}(x)$ in $T_e {\mathbf C}P^n$. As a result we get
a fibre bundle over $G_{k-1} (T_e {\mathbf C}P^n)$ whose fibres are diffeomorphic to
${\mathbf C}P^{n-k}$. Let us denote  this fibre bundle by $T^1 _{k-1,n}$. Obviously,
$T^1 _{k-1,n}$ is also a fibre bundle over the complex Grassmannian
$G_k (T_e {\mathbf C}P^n)$ with the natural projection $p: (v, x) \mapsto v\wedge x $. So we
have the following fibrations
$$ {\mathbf C}P^{k-1} \longrightarrow T^1 _{k-1,n} \longrightarrow G_k (T_e {\mathbf C}P^n) ,$$
$${\mathbf C}P^{n-k} \longrightarrow T^1 _{k-1,n} \longrightarrow G_{k-1} (T_e {\mathbf C}P^n ) .$$
We observe that the invariant metric on $T^1 _{k-1,n} \simeq U_n /(U_{k-1}
\times U_{n-k} \times U_1 )$, obtained from the bi-invariant metric on $U_n$ 
factorized by the action of its subgroup $U_{k-1} \times U_{n-k} \times U_1 $, coincides
with those which are obtained by lifting the invariant metric on $G_{k-1}
(T_e {\mathbf C}P^n)$ via $q$, and the  one on $G_k (T_e {\mathbf C}P^n)$ via $p$. Therefore we
get
$$cd (e, V^{2k} ) = A_{k,n} \int_{G_{k-1}( T_e ({\mathbf C}P^n ))} \int_{{\mathbf C}P^{n-k} (y)}
| < \overline{V^{2k}} , \overline{y \wedge x }> | \, dx \, dy ,  $$
where $A_{k,n}$ is a constant which depends only on $n$ and $k$.

For any point $y \in G_{k-1} (T_e {\mathbf C}P^n)$ denote $\Pi _V y $ the orthogonal
projection of $y$ on the subspace $V^{2k}$. Let $\Pi _V y ^{\perp}$ 
denote the orthogonal complement to the projection $\Pi _V y$ in $V^{2k}$. Then
we get
$$ \int_{{\mathbf C}P^{n-k} (y)} |< \overline{V^{2k}}, \overline{y \wedge x} >| \, dx =
|<V, \overline{y} >| \cdot \int_{{\mathbf C}P^{n-k} (y)} |< \overline{\Pi _V y ^{\perp} }, \overline{x}>| \, dx .
\eqno (3.11)$$
From the proof of Proposition 2.11 we conclude that the right hand side of (3.11)
is less than or equal to $|<V, \overline{y} >|$. Moreover, the equality holds if and
only if $\Pi _V y ^{\perp}$ is a complex line. Repeating the reduction
procedure as above we obtain  Proposition 3.10 from the following lemma.

{\bf Lemma 3.11}. {\it Let $V^{2k}$ be a subspace of real dimension $2k$ in ${\mathbf C}^{
n+1}$.
For every $x \in {\mathbf C}P^n$ let us denote $|< V^{2k},x>|$ the volume of the
projection of the unit complex line $x \in {\mathbf C}P^n$ on the space $V^{2k}$. Then the function
$$ M_C (V^{2k}) = \int_{{\mathbf C}P^n} |<V^{2k} , x>|\, dx $$
reaches its maximal value if and only if $V^{2k}$ is a complex subspace.}

{\bf Proof.} We consider the Hopf fibration $S^{2n+1} \longrightarrow {\mathbf C}P^n$. As in the proof of Proposition 2.11 we conclude that
$$ M_C (V^{2k}) = C_n \int_{S^{2n+1}} |< V^{2k}, x'\wedge Jx'>| \, dx' = C_n \int_{S^{2n+1}} \text{  vol}( \Pi _V x' \wedge 
\Pi _V Jx' ) \, dx' ,$$
where $C_n= \text{  vol} (U_1)^{-1}$, and $\Pi _V x'$ denotes the orthogonal projection of the unit vector $x'\in S^{2n+1}$ on the subspace $V^{2k}$.
Therefore we obtain
$$ M_C (V^{2k}) \le C_n \cdot \int_{S^{2n+1}} |\Pi _V x'| \cdot | \Pi _V Jx' | \, dx ',
\eqno (3.12)$$
and besides, the equality holds iff $\Pi _V x'$ is perpendicular to
$\Pi _V Jx'$ for every $x' \in S^{2n+1}$. That condition is equivalent to
the complexity of $V^{2k}$. Note that the group $SO_{2n+2}$ acts on
the Grassmannian of real $2k$-dimensional planes in ${\mathbf R}^{2n+1} = {\mathbf C}^{n+1}$ transitively. Applying the Schwarz inequality for integrals to the right hand side of (3.12) we get
$$M_C (V^{2k}) \le C_n(\int_{S^{2n+1}} |\Pi _V x'|^2 \, dx')^{1/2}(\int_{S^{2n+1}} |\Pi_V Jx'|^2 \, dx')^{1/2} = C_n \int_{S^{2n+1}} |\Pi _ V x'| ^2 \, dx'  .$$
Moreover, the inequality becomes   an equality if and only if $V$ is a complex
plane (and in this case we also have $|\Pi _V x'| = | \Pi _V Jx'|$). This completes
the proof of Lemma 3.11 and then the proof of Proposition 3.10.

{\it  Continuation of Proof of Theorem 3.2}. The remaining part of this proof
can be carried out in the same way as in the proof of Theorem 3.1. It is easy
to see that the following key lemma  is an analog of Proposition 3.4.

{\bf Lemma 3.12}. {\it For each real plane $V^{2pm} \subset \oplus _{r=0} ^{m-1}
I^r ({\mathbf C}^q)$
 we put
$$M(V) = \int_{{\mathbf C}P^{q-1}} |<\overline{V}, x \wedge \cdots \wedge I^{m-1}(x)>|\, dx .$$ 
Then M(V) reaches its maximal value if and only if $V = V_1 \wedge  \cdots  \wedge I^{m-1} (V_1)$, where   $V_1$ is some complex subspace in ${\mathbf C}^q$.}

{\bf Proof}. Applying the Schwarz inequality and the technique in the proof
of Proposition 3.4 we get
$$M(V) \le C_{q,m} (\int_{S^{2q-1}} B(x,x)^{m/2} \, dx ) (\int_{S^{2q-1}}
B(Jx,Jx)^{m/2} \, dx ) ,$$
where $C_{q,m}$ is some constant and $B(x,x)$ is a symmetric bilinear form as in the proof of Proposition 3.4.
Now, the condition that $M(V^{2k})$ reaches its maximal value is the combination of the following two: $V^{2k}$ is product of $I^r ({\mathbf R}^{2p})$ and $V^{2k}$ is a complex
subspace. This completes the proof of Lemma 3.12.

{\bf Proof of Theorem 3.3}. We follow the proof of Theorem 3.2. 
To do this we consider the Hopf fibration $S^{4q-1} \longrightarrow {\mathbf H}P^{q-1}$
and apply the H\"older inequality for integrals (instead of the Schwarz
inequality).

\vspace{0,5cm}

{\bf \S 4. Properties of $(M)^*$-minimal cycles.}

Let $N$ be a $k$-cycle in Riemannian manifold $M^m$ provided with a family
$(M)^*$ of submanifolds $N^* _{\lambda}$ in $M$ realizing a cycle $[N^*]$ 
as in Corollary 2.2.
If the inequality in this corollary for the volume of $N$ becomes   an
equality, we will call $N$ a $(M)^*$-minimal cycle. Corollary 2.3 states
that a $(M)^*$-minimal cycle is  homologically volume-minimizing. The homological
class $[N] \in H_* (M)$ of such a cycle will be called a $(M)^*$-class. First we
show that there is an analog of {\it Equidistribution Theorem} for homologically
volume-minimizing cycles in a $(M)^*$-homology class.

{\bf Theorem 4.1}. {\it Equidistribution Theorem. Let $N'$ be a homological
volume-minimizing cycle in a $(M)^*$-homology class. Then the set of 
$N^* _{\lambda}
\in (M)^*$ such that $ \# (N^* _{\lambda} \cap N' )\not=  \chi$ is of measure
zero in $(M)^*$. Here $\chi$ equals the intersection number of cycles $[N]$ and $[N^* ]$.}

{\bf Proof}. By our assumption and taking into account Corollary 2.2 we 
conclude that $N'$ also satisfies the condition in Corollary 2.3. Namely we
have
$$ \text{  vol} (N') = \chi \cdot (\sigma (M)^* _k)^{-1} \cdot \text{  vol}(M)^*.$$
Theorem 2.1 implies that the above equality holds if
and only if $N'$ satisfies the following two conditions

1)For almost all $x \in N'$ we have $cd(x, T_x N) = \sigma (M)^*_k$.\\
2)For almost all $y \in (M)^*$ the actual intersection number $\# (N_y \cap
 N')$ equals the algebraic intersection number $\chi$.

Now Theorem 4.1 follows from the second condition.

Applying Theorem 4.1 to
complex submanifolds in the complex projective manifolds  ${\mathbf C}P^n$ we obtain the
following corollary. Recall that the homology group $H_{2k}({\mathbf C}P^n, {\mathbf Z}) = {\mathbf Z}$ is generated by the 
element $[{\mathbf C}P^k]$.

{\bf Corollary 4.2}. {\it Let $r$ be a positive integer, and let $N^{2k}$ be a complex submanifold realizing the   
element $r [{\mathbf C}P^k] \in H_{2k} ({\mathbf C}P^n,{\mathbf Z})$. Then the set of $(2n-2k)$-dimensional
projective spaces ${\mathbf C}P^{n-k} _{\lambda}\subset {\mathbf C}P^n$ such that $\# ({\mathbf C}P^{n-k} _{\lambda}
\cap  N^{2k}) \not= r$ is of measure zero in the set of all ${\mathbf C}P^{n-k}_{\lambda}$ which is identified with $
SU_n /S(U_{n-k} \times U_k)$ provided with the invariant measure}.

{\bf Proof.} Applying Proposition 2.10.a to the cycle $r \, {\mathbf C}P^k$ we get that all
homology classes in $H_* ({\mathbf C}P^n,Z)$ are $(M)^*$-homology classes. It is well known that the complex submanifold $N^{2k}$ is volume-minimizing in its
homology class. Hence we infer Corollary 4.2 from Theorem 4.1.

Volume-minimizing cycles in an $(M)^*$-homology class possess some properties
similar to those of $\phi$-currents, where $\phi$ is a calibration on $M$.
First, we note that the cycles under consideration are also $(M)^*$-minimal.
Further, the tangent space to a $(M)^*$-minimal cycle belongs to a certain
distribution of $k-$planes in $TM$. Namely at every point $x \in M$ we put
$$I(x) = \{ V \in G_k (T_x M) |  \; cd(x,V) = \sigma (M)^* _k \} . $$
Then $(M)^*$-minimal cycles are integral submanifolds of the distribution $I(x)$. Recall that $\phi$-submanifolds are integral submanifolds of the distribution $G_{\phi} (M)$ $= \{V \in TM | \; \phi (\overline{V}) =1 \} $. When $M= G/H$ is
a compact homogeneous Riemannian space, we find a striking relation between
these distributions. Let $\phi $ be an invariant calibration on $M$. Then
its restriction to the tangent space of $M$ at the point $\{ eH\}$ is a $H$-invariant
form. Thererfore, the
value of $\phi$ at a $k$-vector $\overline{V}  \subset T_{\{ eH\} } G/H$ can be expressed as follows
$$ \phi (\overline{V}) = \int_H <\overline{V}, Ad_{\tilde{x}}\overline{W}> \, d\tilde{x}  ,$$
where $\overline{W}$ is some $k$-vector in the space $T_{\{ eH\} }M$. Obviously, the value $\phi (\overline{V})$
depends only on the orbits of the $H$-action on $\bigwedge _k T_{\{ eH\} } M$ (cf. Proposition 2.6). Moreover, let us denote $L$ the isotropy group of the $H$-action at the $k$-vector $\overline{W}$. Then we have
$$\phi (\overline{V}) = \int_{H/L} <\overline{V}, Ad_{\tilde x}\overline{ W}>\, dx  .\eqno(4.1) $$
This formula is similar to the one we used for computing deformation coefficient $cd(\{ eH\} ,V)$,
(see (2.8)).  Further, the distribution $G_{\phi}$ is the set of all $k$-dimensional tangent subspaces whose associated unit simple
k-vectors maximize $\phi (\overline{V})$ ; the distribution $I$ is the set
of all $k$-dimensional tangent subspaces whose associated unit simple $k-$vectors maximize value $cd(x,\overline{V})$. In many cases,
for example, for a K\"ahler form and its powers $\phi$, we can choose
a corresponding $\overline{W}$ as a simple polyvector.  

The similarity between $(M)^*$-cycle and $\phi$-currents also appears in the
following theorem.

{\bf Theorem 4.3.} {\it Let $N$ be a $(M)^*$-minimal cycle realizing a
torsion free
element in the homology group $H_k(M,{\mathbf Z})$. If $M$ is a compact
manifold, then $N$ is a $\phi$-current
for some calibration $\phi$ on $M$ and the homology class $[N]$ is stable.}

{\it Remark}. In many cases, for example, for $M ={\mathbf C}P^n$, there is a unique (up to multiplication by a constant) invariant
calibration of a given dimension on the manifold $M$ (see also [Le 4]).
In such cases, in view of Theorem 4.3, we can obtain a calibration on $M$ with
the help of integral geometry. As it was discussed above, the two kinds of
involved integral inequalities are similar but not equivalent. For instance,
we consider the deformation coefficient as in Proposition 3.5. It is easy
to see that if $m$ is even, then the integrand $|<\overline{V^{km}}, Ad_{\tilde{x}} \overline{W}>|$ equals
$<\overline{V^{km}}, Ad_{\tilde{x}} \overline{W}>$ for all $V^{km}$ which belongs to the distribution of maximal deformation
coefficient. Therefore, such a plane $V^{km}$ also belongs to the distribution of the calibration
associated with $W$ as it was discussed above (see (4.1)).

{\bf Proof of Theorem 4.3}. Let us recall the Federer Stability Theorem. 

{\bf Theorem.} [Fe 2 ]. {\it For every $\alpha \in H_k (M,G)$ we put }
$$\text{mass}(\alpha )= \min \{ \text{vol} X^k \subset M \vert \; [X^k] = \alpha \} .$$ 
{\it Then the following
equality holds for $\alpha \in H_k (M, {\mathbf Z})$.}

$$\lim_{n \to \infty} { \text{ mass}(n \alpha) \over n } = \text{mass}(\alpha _{\mathbf R}) ,$$

{\it where $\alpha _{\mathbf R}$ denotes the image of $\alpha$ under the map $H_k (M,{\mathbf Z}) 
\rightarrow H_k (M,{\mathbf R})$. }

If for some $n \in {\mathbf Z}^+$ we have $ \text{  mass} (n \alpha) /n = \text{  mass}(\alpha _{\mathbf R})$ we say
that the homology class $\alpha$ is stable. 

Now assume $N$ is as in Theorem 4.3. We observe that the cycle $pN$ is 
also a $(M)^*$-cycle for all $p \in {\mathbf Z}^+$. So we get
$$\text{  mass}(  p[N] ) /p = \text{  mass} ([N]).$$
Therefore, according to the Federer Stability Theorem, the homology class $[N]$ must be stable, and $N$ is a
volume-minimizing cycle in the class $[N]_{\mathbf R} \in H(M,{\mathbf R})$. It is well-known that there is a calibration $\phi$ on $M$ which calibrates $N$ (cf. [D-F], [Le 4]).

Applying Theorem 4.3 to Theorem 3.1 we obtain the following corollary.

{\bf Corollary 4.4}.[G-M-Z]. {\it If the Grassmannian of oriented planes $G^+ _k
({\mathbf R}^{k+m})$ realizes a non-trivial element in the homology group $H_{km} (
G^+ _l ({\mathbf R}^{l+m}), {\mathbf R})$ with real coefficients, then $G^+ _k ({\mathbf R}^{k+m})$
is a volume-minimizing cycle in its homology class with real coefficients.}  

{\bf Proof}. Obviously, $G _k( {\mathbf R}^{k+m})$ and its 2-sheeted covering $G^+ _k
({\mathbf R}^{k+m})$ have the same homology groups with real coefficients. By Theorem 4.3,
$G_k ({\mathbf R}^{k+m})$ is a volume-minimizing real current. Its is well known that
in this case there exists an invariants calibration $\phi$ on $G_l( {\mathbf R}^{l+m})$
such that $\phi$ calibrates $G_k( {\mathbf R}^{k+m})$. It is easy to see that the
lifted calibration $\phi ^*$ on $G^+ _l ({\mathbf R}^{l+m})$ must calibrate $G^+ _k 
({\mathbf R}^{k+m})$ too. This means that $G^+ _k ({\mathbf R}^{k+m})$ is a globally minimal
submanifold.

Finally we conjecture that every homology class in $H_8 (F_4/Spin_9, {\mathbf Z})$
is a $(M)^*$-class. A. T. Fomenko and M. Berger proved that the Helgason sphere $S^8$ realizing
the generating element of this group is a globally minimal submanifold [Fo 1],
[Be] . We also conjecture that every canonically embedded sub-Grassmannian
$G_k({\mathbf F}^l) \subset G_{k+m}({\mathbf F}^{l+n})$ is volume minimizing in its ${\mathbf Z_2}$ homology,
where ${\mathbf F}= \mathbf {R,C,H}$ (see also [G-M-Z] for the case of oriented $G^+ _{k+m}
({\mathbf R}^{l+n})$).

\vspace{0,5cm}

{\bf Acknowledgement}. This paper was conceived during my visit to the
International Centre for Theoretical Physics at Trieste and now has been 
completed 
in the Max-Planck Institut f\"ur Mathematik at Bonn. I would like to thank the
ICTP and MPI f\"ur Mathematik for hospitality and financial support. I am grateful
to Prof. H. Karcher for discussion and pointing out the nice book [Sa],
and the referee for critical comments and useful suggestions.

\

{\bf Appendix:  Correspondence  with Professor  Tasaki  on Proposition 3.5}

From tasaki\@ math.tsukuba.ac.jp Fri Sep  2 08:34:41 1994\\
To: lehong\@ mpim-bonn.mpg.de\\
Subject: Question

Dear Professor Le,

   I have been reading your paper "Application of integral geometry 
to minimal surfaces" with great interest.  I gave a lecture on integral
geometry which included your results in the paper and mine.  At that
time there was a point which I did not understand.  I would like to
continue to give such a lecture, so I hope to make it clear.

   In the proof of Proposition 3.5, I think, you do not prove that 
$V_0 = Ad_xW$ if the deformation coefficient $cd(e, V)$ attains its
maximum at $V_0$ in the case of $m = 2$.  In this case we can define a
complex structure I' which coincides with I on $R_1^l$ and $-I$ on $ R_2^l$.
Then the equalities of (3.4) and (3.5) hold if $V$ is a complex subspace
of complex dimension 1, which may not be of the form $Ad_xW$.

   I am looking forward to hearing from you.

   Sicerely yours,
   
   Hiroyuki Tasaki

\

From lehong Tue Sep  6 10:37:09 1994\\
To: tasaki\@ math.tsukuba.ac.jp\\
Subject: Re: Question

Dear Professor Tasaki,

Thank you very much for your mail. Certainly I overlooked the case
$m=2$. But it is not hard to correct the classification theorem 3.1' since
there is a natural Hermit structure on $G_k (R^{m+k})$ which you already noticed.
(by the way G-M-Z also classified  for $m\ge 4$). The correct statement 
should be so: if $m=2$ then $M$ must be a Hermit submanifold.
Proof: Clearly the class $[G_k (R^{m+k})]$ is a $(M^*)$-class. By Theorem 4.3
this class is stable, in particular $M$ is a minimizing real current. Since
the Wirtinger form $\Omega$ calibrates $G_k (R^{m+k})$ this $\Omega$ also calibrates
$M$. Hence $M$ is a complex submanifold.

Remark 1. The proof goes through for both coefficient group $\mathbf Z$ and $\mathbf Z_2$.

Remark 2. I suspect that the  Proposition 3.5 (corrected for the case
$m=2$) should include all complex planes (of dimension $p$) but have yet not
proper proof (it is easy to see that complex planes satisfy the condition
but the other side is more complicated. )

Best regards,

Le Hong Van.

\

From tasaki\@ math.tsukuba.ac.jp Mon Sep 12 13:41:32 1994\\
To: lehong\@ mpim-bonn.mpg.de\\

Dear Professor Le,

Thank you for your clear reply.  When I give a lecture on integral geometry,
can I use the result mentioned in your e-mail?

Do you know DIFFERENTIAL GEOMETRY E-PRINTS of MSRI?  It may be usefull for us.
If you send an e-mail to dg-ga\@ msri.org only with Subject: help, then you can 
get information about it.

Sincerely yours,

Hiroyuki Tasaki

\

From lehong Wed Sep 14 11:02:26 1994\\
To: tasaki\@ math.tsukuba.ac.jp\\
Subject: Re: Question

Dear Professor Tasaki,

Thank you very much for your mail. Certainly I would be very happy
if you include my result in your lecture. I will try Diff. Geom.
E-Prints of MSRI.
I hope we will meet again in the future.

Best wishes,

Le Hong Van.










\vspace{1cm}
\begin{center}
{\sc REFERENCES}
\end{center}
\vspace{1cm}

[Be] M. Berger, {\it Du cot\'e de chez Pu}, Ann. Scient. Ec. Norm. Sup.
4,(1972), 1-44.

[Ch] S. S. Chern, {\it Complex Manifolds without Potential Theory,} 
Springer-Verlag, NewYork-Heidenberg-Berlin, 1979.

[D] Dao Chong Thi, {\it Minimal real currents on compact Riemannian 
manifolds}, Math. USSR Izvest., 11 (1977), 807-820

[D-F] Dao Chong Thi and A. T. Fomenko, {\it Minimal Surfaces and the
Plateau Problem}, A.M.S., Providence, R.I., 1991.

[DGGW 1] D. DeTruck, H. Gluck, C. Gordon, D. Webb,
{\it The geometry of isospectral deformation}, preprint (1990).

[DGGW 2] D .DeTruck, H. Gluck, C. Gordon, D. Webb,
{\it  Inaudible geometry of nilmanifolds}, preprint (1990).

[Fe 1] H. Federer, {\it Geometric Measure Theory}, Springer-Berlin, 1969.

[Fe 2] H. Federer, {\it Real flat chains, cochains and variational problems}, 
Ind. Univ. Math. J., 26 (1974), 351-467.

[Fo 1] A. T. Fomenko, {\it Minimal compacta in Riemanian manifolds and 
Reifenberg's conjecture}, Math. USSR Isv. 6(1972), 1037-1066.

[Fo 2] A. T. Fomenko, {\it Variational Principles in Topology}, Kluwer
Academic Publisher, Dordrecht, 1990.

[G-M-Z] H. Gluck, F. Morgan, W. Ziller, {\it Calibrated geometries in
Grassmann manifolds}, Comm. Math. Helv. 64 (1989), 256-268.

[G-M-M] H. Gluck, D. Mackenzie and F. Morgan, {\it Volume minimizing
cycles in Grassmann manifolds} , preprints 1991.

[H-L] R. Harvey and H. B. Lawson, {\it Calibrated geometries}, Acta Math.,
148 (1982), 47-157.

[Hs-Ln] W-Y. Hsiang and H. B. Lawson, {\it Minimal submanifolds of low
cohomogeneity}, J. Diff. Geom. 5 (1971), 1-38.

[K-R] V. E. Katsnelson and L. I. Ronkin, {\it On the minimal volume of an analytic set},
Sibirs. Math. J., 15 (1974), 370-378.

[Ln 1] H. B. Lawson, {\it The equivariant Plateau problem and interior
regularity}, Trans. Am. Math. Soc., 173 (1973), 231-249.

[Ln 2] H. B. Lawson, {\it The stable homology of a flat torus,} Math. Scand.
36 (1972), 49-73.h. 95 (1989), 437-446.

[Le 1] Le Hong Van, {\it Relative calibrations and stability of minimal 
surfaces}, Lect. Notes in Math., Springer-Verlag, v.1453, (1990), 245-262.  

[Le 2] Le Hong Van, {\it Minimal $\Phi$-Lagrangian surfaces in 
almost Hermitian manifolds}, Math. USSR Sbornik, v.67, (1990), 379-391. 

[Le 3] Le Hong Van, {\it Application of integral geometry to minimal
surfaces}, preprint MPI/92-7.

[Le 4] Le Hong Van, {\it Each globally minimal surfaces in a compact
homogeneous space has an invariant calibration}, Soviet Math. Dokl., 41,
(1990), 66-69.

[Le-Fo] Le Hong Van and A. T. Fomenko {\it Volumes of minimal surfaces and curvature tensor of Riemannian manifolds}, Soviet Math. Dokl., 37, (1988),
817-820.

[Sa] L. A. Santal\'o, {\it Integral Geometry and Geometric Probability},
Addison-Wesley Publishing Company, U.S.A., 1976.

[Si] J. Simons, {\it Minimal varieties in Riemannian manifolds},
Ann. Math., 88, (1968), 62-105.

[T] H. Tasaki, {\it Calibrated geometries in quaternionic Grassmannian},
Osaka J. Math. 25 (1988), 591-597.

\vspace{1cm}





\end{document}